\input amstex
\documentstyle{amsppt}
%
\catcode`@=11
\redefine\output@{%
  \def\break{\penalty-\@M}\let\par\endgraf
  \ifodd\pageno\global\hoffset=105pt\else\global\hoffset=8pt\fi  
  \shipout\vbox{%
    \ifplain@
      \let\makeheadline\relax \let\makefootline\relax
    \else
      \iffirstpage@ \global\firstpage@false
        \let\rightheadline\frheadline
        \let\leftheadline\flheadline
      \else
        \ifrunheads@ 
        \else \let\makeheadline\relax
        \fi
      \fi
    \fi
    \makeheadline \pagebody \makefootline}%
  \advancepageno \ifnum\outputpenalty>-\@MM\else\dosupereject\fi
}
\def\Beta{\mathchar"0\hexnumber@\rmfam 42}
\catcode`\@=\active
\nopagenumbers
\def\negskp{\hskip -2pt}

\accentedsymbol\hatgamma{\kern 2pt\hat{\kern -2pt\gamma}}
\accentedsymbol\checkgamma{\kern 2.5pt\check{\kern -2.5pt\gamma}}
\def\blue#1{#1}

\catcode`#=11\def\diez{#}\catcode`#=6
\catcode`&=11\catcode`&=4
\catcode`_=11\def\podcherkivanie{_}\catcode`_=8
\def\mycite#1{\cite{\blue{#1}}\immediate\special{ps:
     ShrHPSdict begin /ShrBORDERthickness 0 def}}
\def\myciterange#1#2#3#4{\cite{\blue{#2#3#4}}\immediate\special{ps:
     ShrHPSdict begin /ShrBORDERthickness 0 def}}
\def\mytag#1{%
    \tag#1}
\def\mythetag#1{\thetag{\blue{#1}}\immediate\special{ps:
     ShrHPSdict begin /ShrBORDERthickness 0 def}}
\def\myrefno#1{\no#1}
\def\myhref#1#2{\blue{#2}\immediate\special{ps:
     ShrHPSdict begin /ShrBORDERthickness 0 def}}
\def\myEarXivlink{\myhref{http://arXiv.org}{http:/\negskp/arXiv.org}}

\def\mytheorem#1{\csname proclaim\endcsname{Theorem #1}}
\def\mytheoremwithtitle#1#2{\csname proclaim\endcsname{Theorem #1#2}}
\def\mythetheorem#1{\blue{#1}\immediate\special{ps:
     ShrHPSdict begin /ShrBORDERthickness 0 def}}
\def\mylemma#1{\csname proclaim\endcsname{Lemma #1}}
\def\mylemmawithtitle#1#2{\csname proclaim\endcsname{Lemma #1#2}}
\def\mythelemma#1{\blue{#1}\immediate\special{ps:
     ShrHPSdict begin /ShrBORDERthickness 0 def}}
\def\mycorollary#1{\csname proclaim\endcsname{Corollary #1}}

\def\myconjecture#1{\csname proclaim\endcsname{Conjecture #1}}
\def\myconjecturewithtitle#1#2{\csname proclaim\endcsname{Conjecture #1#2}}
\def\mytheconjecture#1{\blue{#1}\immediate\special{ps:
     ShrHPSdict begin /ShrBORDERthickness 0 def}}

\pagewidth{360pt}
\pageheight{606pt}
\topmatter
\title
A note on the first cuboid conjecture.
\endtitle
\author
Ruslan Sharipov
\endauthor
\address Bashkir State University, 32 Zaki Validi street, 450074 Ufa, Russia
\endaddress
\email\myhref{mailto:r-sharipov\@mail.ru}{r-sharipov\@mail.ru}
\endemail
\abstract
    Recently the problem of constructing a perfect Euler cuboid was related 
with three conjectures asserting the irreducibility of some certain three 
polynomials depending on integer parameters. In this paper a partial result
toward proving the first cuboid conjecture is obtained. The polynomial which,
according to this conjecture, should be irreducible over integers is proved 
to have no integer roots. 
\endabstract
\subjclassyear{2000}
\subjclass 11D41, 11D72, 12E05\endsubjclass
\endtopmatter
\TagsOnRight
\document

\head
1. Introduction.
\endhead
     An Euler cuboid is a rectangular parallelepiped whose edges and 
face diagonals all are of integer lengths. A perfect cuboid is an Euler 
cuboid whose space diagonal is also of an integer length. Cuboids with 
integer edges and face diagonals are known since 1719 (see 
\myciterange{1}{1}{--}{35}), however, no perfect cuboid is known by now. 
The problem of constructing perfect cuboids or proving their non-existence 
is an open mathematical problem.\par
     In \mycite{36} the problem of constructing perfect cuboids was reduced
to the polynomial Diophantine equation $P_{abu}(t)=0$, where $P_{abu}(t)$ is
given by the formula
$$
\gathered
P_{abu}(t)=t^{12}+(6\,u^2\,-2\,a^2\,-2\,b^2)\,t^{10}
+(u^4\,+b^4+a^4+4\,a^2\,u^2+\\
+\,4\,b^2\,u^2-12\,b^2\,a^2)\,t^8+(6\,a^4\,u^2+6\,u^2\,b^4-8\,a^2\,b^2\,u^2-\\
-\,2\,u^4\,a^2-2\,u^4\,b^2-2\,a^4\,b^2-2\,b^4\,a^2)\,t^6+(4\,u^2\,b^4\,a^2+\\
+\,4\,a^4\,u^2\,b^2-12\,u^4\,a^2\,b^2+u^4\,a^4+u^4\,b^4+a^4\,b^4)\,t^4+\\
+\,(6\,a^4\,u^2\,b^4-2\,u^4\,a^4\,b^2-2\,u^4\,a^2
\,b^4)\,t^2+u^4\,a^4\,b^4.
\endgathered\quad
\mytag{1.1}
$$
The main result of \mycite{36} is formulated in the following theorem.
\mytheorem{1.1} A perfect Euler cuboid does exist if and only if the Diophantine 
equation $P_{abu}(t)=0$ has a solution such that $a$, $b$, $u$, and $t$ are 
positive integer numbers obeying the inequalities $t>a$, $t>b$, $t>u$, and 
$(a+t)\,(b+t)>2\,t^2$\!.
\endproclaim
     Note that $P_{abu}(t)$ is a polynomial of four variables
$a$, $b$, $u$ and $t$. However, in the formula \mythetag{1.1} it is presented as 
a univariate polynomial depending on three integer parameters $a$, $b$, and $u$. 
Relying on this presentation, in \mycite{37} the theorem~\mythetheorem{1.1} was 
reformulated as follows. 
\mytheorem{1.2} A perfect Euler cuboid does exist if and only if for some
positive coprime integer numbers $a$, $b$, and $u$ the polynomial equation 
$P_{abu}(t)=0$ has a ra\-tional solution $t$ obeying the inequalities \pagebreak 
$t>a$, $t>b$, $t>u$, and $(a+t)\,(b+t)>2\,t^2$\!. 
\endproclaim
     If the equation $P_{abu}(t)=0$ has a rational solution, then the polynomial
\mythetag{1.1} with integer coefficients is reducible over the field of rational 
numbers. Note that the leading coefficient of this polynomial is equal to unity. 
Hence due to the rational root theorem (see \mycite{38}, \mycite{39}, or 
\mycite{40}) each rational root of the polynomial $P_{abu}(t)$, if any, is 
necessarily integer and $P_{abu}(t)$ is reducible over the ring of integers.\par
     In \mycite{37} the polynomial \mythetag{1.1} was studied for reducibility 
and the following special cases were discovered where $P_{abu}(t)$ is reducible:
$$
\xalignat 3
&\hskip -2em
\text{1) \ }a=b; &&\text{3) \ }b\,u=a^2; &&\text{5) \ }a=u;
\qquad\\
\vspace{-1.5ex}
\mytag{1.2}\\
\vspace{-1.5ex}
&\hskip -2em
\text{2) \ }a=b=u; &&\text{4) \ }a\,u=b^2; &&\text{6) \ }b=u.
\qquad
\endxalignat
$$
Being reducible in the cases \mythetag{1.2}, the polynomial \mythetag{1.1}
gives rise to the polynomials 
$$
\gather
\hskip -2em
\gathered
P_{au}(t)=t^8+6\,(u^2-a^2)\,t^6+(a^4-4\,a^2\,u^2+u^4)\,t^4-\\
-\,6\,a^2\,u^2\,(u^2-a^2)\,t^2+u^4\,a^4,
\endgathered
\mytag{1.3}\\
\vspace{1ex}
\hskip -2em
\gathered
Q_{pq}(t)=t^{10}+(2\,q^2+p^2)\,(3\,q^2-2\,p^2)\,t^8
+(q^8+10\,p^2\,q^6+\\
+\,4\,p^4\,q^4-14\,p^6\,q^2+p^8)\,t^6
-p^2\,q^2\,(q^8-14\,p^2\,q^6+4\,p^4\,q^4+\\
+\,10\,p^6\,q^2+p^8)\,t^4-p^6\,q^6\,(q^2+2\,p^2)
\,(-2\,q^2+3\,p^2)\,t^2-q^{10}\,p^{10}
\endgathered
\mytag{1.4}
\endgather
$$
depending on the integer parameters $a$, $u$ and $p$, $q$. In \mycite{37}
the reducibility of the polynomials \mythetag{1.3}, \mythetag{1.4} and 
the reducibility of the initial polynomial \mythetag{1.1} were studied 
numerically and three conjectures were formulated. 
\myconjecture{1.1} For any positive coprime integers $a\neq u$ the polynomial
$P_{au}(t)$ in \mythetag{1.3} is irreducible in the ring $\Bbb Z[t]$. 
\endproclaim
\myconjecture{1.2} For any positive coprime integers $p\neq q$ the 
polynomial $Q_{pq}(t)$ in \mythetag{1.4} is irreducible in the ring 
$\Bbb Z[t]$. 
\endproclaim
\myconjecture{1.3} For any three positive coprime integer numbers $a$, $b$, 
and $u$ such that none of the conditions \mythetag{1.2} is satisfied the 
polynomial $P_{abu}(t)$ in \mythetag{1.1} is irreducible in the ring $\Bbb Z[t]$. 
\endproclaim
    The main goal of this paper is to prove the following partial result
associated with the first cuboid conjecture~\mytheconjecture{1.1}.
\mytheorem{1.3} For any positive coprime integers $a\neq u$ the polynomial
$P_{au}(t)$ in \mythetag{1.3} has no integer roots. 
\endproclaim
\head
2. The inversion symmetry and parity.
\endhead
     The polynomial $P_{au}(t)$ in \mythetag{1.3} possesses some special property.
It is expressed by the following formula which can be verified by direct calculations:
$$
\hskip -2em
P_{au}(t)=\frac{P_{au}(i\,a\,u/t)\,t^8}{a^4\,u^4}.
\mytag{2.1}
$$
Here $i=\sqrt{-1}$. The formula \mythetag{2.1} contains the inversion of $t$ in 
$P_{au}(t)$. \pagebreak For this reason I call it the inversion symmetry. Apart 
from \mythetag{2.1}, we have 
$$
\hskip -2em
P_{au}(t)=P_{au}(-t).
\mytag{2.2}
$$
The formula \mythetag{2.2} means that the polynomial $P_{au}(t)$ is even. 
\head
3. Breaking the proof of irreducibility into special cases. 
\endhead
      The irreducibility of polynomials is usually proved by contradiction.
If the conjecture~\mytheconjecture{1.1} is not valid, this would mean that the
polynomial \mythetag{1.3} is reducible, i\.\,e\. it is presented as a product
of two non-constant polynomials
$$
\hskip -2em
P_{au}(t)=A(t)\,B(t).
\mytag{3.1}
$$
Since $\deg P_{au}(t)=8$, the equality \mythetag{3.1} assumes four special cases:
$$
\xalignat 2
&\hskip -2em
\text{1) \ }P_{au}(t)=A_1(t)\,B_7(t),
&&\text{2) \ }P_{au}(t)=A_2(t)\,B_6(t),\\
\vspace{-1.5ex}
\mytag{3.2}\\
\vspace{-1.5ex}
&\hskip -2em
\text{3) \ }P_{au}(t)=A_3(t)\,B_5(t),
&&\text{4) \ }P_{au}(t)=A_4(t)\,B_4(t).
\endxalignat
$$
Other three cases $P_{au}(t)=A_5(t)\,B_3(t)$, $P_{au}(t)=A_6(t)\,B_2(t)$, 
$P_{au}(t)=A_7(t)\,B_1(t)$ are equivalent to the cases 1, 2, and 3 up to
the transposition of factors. 
\head
4. The case of a linear factor.
\endhead
     This case is number one in \mythetag{3.2}. In this case $P_{au}(t)=A_1(t)
\,B_7(t)$, where $A_1(t)$ is a linear factor and $B_7(t)$ is its complementary
seventh order factor: 
$$
\hskip -2em
A_1(t)=t-A_0.
\mytag{4.1}
$$
The formula \mythetag{3.1} means that $t=A_0$ is a real integer root of the
polynomial $P_{au}(t)$. Since $a\,\neq 0$ and $b\,\neq 0$, we have $A_0\neq 0$. 
Due to \mythetag{2.1} and \mythetag{2.2}, along with $t=A_0$, the polynomial 
$P_{au}(t)$ has the following real and imaginary roots:
$$
\xalignat 3
&\hskip -2em
t=\frac{i\,a\,b}{A_0}, 
&&t=-A_0, 
&&t=-\frac{i\,a\,b}{A_0}. 
\mytag{4.2}
\endxalignat
$$
The formulas \mythetag{4.1} and \mythetag{4.2} mean that
$$
\hskip -2em
P_{au}(t)=(t^2-A_0^2)\biggl(t^2+\frac{a^2\,u^2}{A_0^2}\biggr)B_4(t). 
\mytag{4.3}
$$
Applying the Gauss's lemma (see \mycite{38}, \mycite{39}, and \mycite{41}), 
we conclude that the fraction $a^2\,u^2/A_0^2$ in \mythetag{4.3} simplifies
to an integer number. Let's denote 
$$
\hskip -2em
C_0=\frac{a\,u}{A_0}.
\mytag{4.4}
$$
Then the formula \mythetag{4.3} is written as follows:
$$
\hskip -2em
P_{au}(t)=\bigl(t^4+(C_0^2-A_0^2)\,t^2-a^2\,u^2\bigr)\,B_4(t)
\text{,\ \ where \ }A_0\,C_0=a\,u. 
\mytag{4.5}
$$\par 
     Now let's apply the formulas \mythetag{2.1} and \mythetag{2.2} to 
\mythetag{4.5}. As a result we get the following symmetries for the
polynomial $B_4(t)$ in \mythetag{4.3} and \mythetag{4.5}:
$$
\xalignat 2
&\hskip -2em
B_4(t)=-\frac{B_4(i\,a\,u/t)\,t^4}{a^2\,u^2},
&&B_4(-t)=B_4(t).
\mytag{4.6}
\endxalignat
$$
The symmetries \mythetag{4.6} mean that the polynomial $B_4(t)$ is given by
the formula
$$
\hskip -2em
B_4(t)=t^4+B_2\,t^2-a^2\,u^2. 
\mytag{4.7}
$$
Substituting \mythetag{4.7} into the formula \mythetag{4.5}, we derive
$$
\hskip -2em
\gathered
P_{au}(t)=t^8+(B_2+C_0^2-A_0^2)\,t^6
+\bigl((C_0^2-A_0^2)\,B_2-2\,a^2\,u^2\bigr)\,t^4-\\
-\,a^2\,u^2\,(B_2+C_0^2-A_0^2)\,t^2+u^4\,a^4.
\endgathered
\mytag{4.8}
$$
Comparing \mythetag{4.8} with the initial formula \mythetag{1.3}, we find 
that
$$
\hskip -2em
\aligned
&B_2+C_0^2-A_0^2=6\,(u^2-a^2),\\
&(C_0^2-A_0^2)\,B_2=(u^2-a^2)^2.
\endaligned
\mytag{4.9}
$$
The equations \mythetag{4.9} should be complemented with the equation
$$
\hskip -2em
A_0\,C_0=a\,u.
\mytag{4.10}
$$
The equation \mythetag{4.10} is taken from \mythetag{4.5}. It is equivalent 
to \mythetag{4.4}. The results of the above calculations are summarized in 
the following lemma.
\mylemma{4.1} For $a\neq 0$ and $u\neq 0$ the polynomial $P_{au}(t)$ in 
\mythetag{1.3} has a linear factor of the form \mythetag{4.1} in the ring 
of polynomials $\Bbb Z[t]$ if and only if the system of Diophantine 
equations \mythetag{4.9} and \mythetag{4.10} is solvable with respect 
to the integer variables $A_0$, $B_2$, and $C_0$. 
\endproclaim
     The Diophantine equations \mythetag{4.9} and \mythetag{4.10} are easily
solvable for $u=\pm\,a$. Indeed, in this case we have the following solution 
for them:
$$
\xalignat 3
&\hskip -2em
A_0=\pm\,a, &&C_0=\pm\,a, &&B_2=0.
\endxalignat
$$
\mylemma{4.2} For $u\neq\pm\,a$ the system of Diophantine equations \mythetag{4.9} 
and \mythetag{4.10} is not solvable with respect to the integer variables $A_0$, 
$B_2$, and $C_0$.  
\endproclaim
\demo{Proof} Let's square the first equation \mythetag{4.9} and let's multiply 
by $36$ the second equation \mythetag{4.9}. As a result we get the equations
$$
\hskip -2em
\gathered
B_2^2+2\,(C_0^2-A_0^2)\,B_2+(C_0^2-A_0^2)^2=36\,(u^2-a^2)^2,\\
36\,(C_0^2-A_0^2)\,B_2=36\,(u^2-a^2)^2.
\endgathered
\mytag{4.11}
$$
Subtracting the second equation \mythetag{4.11} from the first one, we derive
$$
\hskip -2em
B_2^2-34\,(C_0^2-A_0^2)\,B_2+(C_0^2-A_0^2)^2=0.
\mytag{4.12}
$$
The left hand side of the equation \mythetag{4.12} is factored into the product
of two linear terms with respect to $B_2$. As a result this equation is written as
$$
(B_2-(17+12\,\sqrt{2})\,(C_0^2-A_0^2))\,(B_2-(17-12\,\sqrt{2})\,(C_0^2-A_0^2))=0.
\quad
\mytag{4.13}
$$
The equation \mythetag{4.13} breaks into two separate equations, i\.\,e\. 
it means that $A_0$, $B_2$, and $C_0$ should obey one of the following two 
equations:
$$
\hskip -2em
\aligned
&B_2=(17+12\,\sqrt{2})\,(C_0^2-A_0^2),\\
&B_2=(17-12\,\sqrt{2})\,(C_0^2-A_0^2).
\endaligned
\mytag{4.14}
$$
None of the equations \mythetag{4.14} can be satisfied by integer numbers $A_0$, 
$B_2$, and $C_0$ unless $C_0^2=A_0^2$. But if $C_0^2=A_0^2$, from the second
equation \mythetag{4.9} we easily derive $u^2=a^2$ and $u=\pm\,a$. The proof
of the lemma~\mythelemma{4.2} is over. 
\qed\enddemo
Combining the lemmas~\mythelemma{4.1} and \mythelemma{4.2} one easily proves the
following theorem.
\mytheorem{4.1} For any two positive integers $a\neq u$ the polynomial $P_{au}(t)$ 
in \mythetag{1.3} has no linear factors of the form \mythetag{4.1} in the ring 
$\Bbb Z[t]$. 
\endproclaim 
The theorem~\mythetheorem{4.1} implies the theorem~\mythetheorem{1.3} declared 
in the introduction. The theorem~\mythetheorem{1.3} is weaker than the 
conjecture~\mytheconjecture{1.1}. However, if similar results for the other two conjectures~\mytheconjecture{1.2} and \mytheconjecture{1.3} will be obtained, 
this would be sufficient to prove the non-existence of perfect cuboids.

\enddocument
\end